\documentclass{article}
\usepackage[utf8]{inputenc}
\usepackage[english]{babel}
\usepackage{amsmath}
\usepackage{amsthm}
\usepackage{amsfonts}
\usepackage{amssymb}

\begin{document}

\author{Walter Wyss}
\title{On Rational Points on the Elliptic Curve \newline
$E(q) : p^2+q^2 = r^2(1+p^2q^2)$}
\date{}
\maketitle

\begin{abstract}
We look at the elliptic curve $E(q)$, where $q$ is a fixed rational number. A point $(p,r) \in E(q)$ is called a rational point if both $p$ and $r$ are rational numbers. We introduce the concept of conjugate points and show that not both can be rational points.
\end{abstract}

\section{Introduction}
The elliptic curve $E(q)$ has been studied by Edwards\cite{1} together with its Abelian group. This elliptic curve shows up in the study of rational cuboids. A cuboid is a rectangular parallelepiped. It is called perfect if all of its sides, face diagonals, and space diagonals are integers, or equivalently  are rational numbers. A representation of the sides, face diagonals, and space diagonal by trigonometric functions leads to the elliptic curve and the concept of conjugate points on this curve. That conjugate points can not both be rational is due to the fact that there is no perfect cuboid \cite{2}.

\section{Cuboids}
Let $x,y,z$ denote the three different edges. The face rectangle $(x,y)$ has diagonal $a$, the face rectangle $(x,z)$ has diagonal $b$, and the face rectangle $(y,z)$ has diagonal $c$. The space diagonal is denoted by $d$. These quantities satisfy the equations

\begin{align}
 x^2 + y^2 = a^2 \\
 x^2 + z^2 = b^2 \\
 y^2 + z^2 = c^2 \\
 x^2 + y^2 + z^2 = d^2 
\end{align} 
 
These quantities are only determined up to their absolute values. We now introduce two concepts.

\subsection*{\underline{Definition 1}}
The generator of an arbitrary angle $\alpha$ is defined by

\begin{align}
m(\alpha) = \frac{\sin{\alpha}}{1+\cos{\alpha}} = \tan({\frac{\alpha}{2}})
\end{align}

\subsection*{\underline{Definition 2}}
An angle $\alpha$ is called a Heron angle if both $\sin{\alpha}$  and $\cos{\alpha}$ are rational numbers.

\subsection*{\underline{Definition 3}}
The involution of an angle $\alpha$ is given by 
\begin{equation}
\bar{\alpha} = \frac{\pi}{2} - \alpha, \hspace{0.25cm} \bar{\bar{\alpha}} = \alpha
\end{equation}

\subsection*{\underline{Lemma 1}}
\begin{enumerate}
\item Let $m(\alpha)$ be the generator of the angle $\alpha$, then \\
\begin{equation}
\cos{\alpha} = \frac{1-m^2(\alpha)}{1+m^2(\alpha)}, \hspace{0.25cm }\sin{\alpha} = \frac{2m(\alpha)}{1+m^2(\alpha)]}
\end{equation}
\item If $\alpha$ is a Heron angle then $m(\alpha)$ is rational and conversely, if $m(\alpha)$ is rational then $\alpha$ is a Heron angle
\item \begin{equation}
m(\bar{\alpha}) = \frac{1-m(\alpha)}{1+m(\alpha)}, \hspace{0.25cm}
\cos{\bar{\alpha}} = \sin{\alpha}, \hspace{0.25cm }\sin{\bar{\alpha}} = \cos{\alpha} \\
\end{equation}
\item \begin{equation}
m(\pi - \alpha) = \frac{1}{m(\alpha)}
\end{equation}
\end{enumerate}

\begin{proof}
Left to the reader.
\end{proof}

\bigskip

Now let $x,y,z,a,d$ be rational numbers and $\vartheta$ and $\varphi$ be Heron angles. Then we have the representation \\
\begin{align}
x = d \cos{\vartheta}\sin{\varphi} \\
y = d \cos{\vartheta}\cos{\varphi} \\
z = d \sin{\vartheta} \\
a = d \cos{\vartheta} 
\end{align}

Introduce the two angles $\phi, \psi$ through
\begin{align}
\cos{\phi} = \cos{\vartheta}\cos{\varphi} \\
\cos{\psi} = \cos{\vartheta}\sin{\varphi} \\
\text{or} \ \cos{\psi} = \cos{\vartheta}\cos{\bar{\varphi}} 
\end{align}

Then $\cos{\phi}$ and $\cos{\psi}$ are both rational and \\
\begin{align}
x = d\cos{\psi} \\
y = d\cos{\phi} \\
z = d\sin{\vartheta} \\
a = d\cos{\vartheta} \\
b = d\sin{\phi} \\
c = d\sin{\psi} 
\end{align}

Since there is no perfect cuboid, we see that if $b$ is rational ($\sin{\phi}$ is rational) then $c$ can not be rational ($\sin{\psi}$ is not rational) and vice versa. That means that not both $\phi$ and $\psi$ can be Heron angles. \\

We now introduce the generators \\
\begin{equation}
p = m(\varphi), q = m(\vartheta), r = m(\phi) \\
\end{equation}

From equation (14) we find \\
\begin{equation}
\frac{1-r^2}{1+r^2} = \frac{1-q^2}{1+q^2} \frac{1-p^2}{1+p^2}
\end{equation}
or
\begin{equation}
r^2 = \frac{p^2+q^2}{1+p^2q^2}
\end{equation}

For a fixed rational $q$ this is the elliptic curve $E(q)$.

\subsection*{\underline{Definition 4}}
Let $(p,r) \in E(q)$. Then $(\bar{p},s) \in E(q)$ is called its conjugate point.

\section*{Theorem}
If $(p,r) \in E(q)$ is a rational point then its conjugate point $(\bar{p},s) \in E(q)$ is not a rational point and vice versa.

\begin{proof}
If $(p,r) \in E(q)$ is a rational point, then $\phi$ is a Heron angle and $b$ is rational. \\
Since $\bar{p} = \frac{1-p}{1+p}$ is the generator of $\bar{\varphi}$, we see from equation (16) and the fact that $c$ can not be rational, that $\psi$ can not be a Heron angle. This means that $(\bar{p},s) \in E(q)$ is not a rational point. The converse is also true because of the involution.
\end{proof}

\section{The Abelian group belongs to E(q)}
From equation (25) we get \\
\begin{equation}
r^2-p^2 = q^2 (1-r^2p^2) , (p,r)\in E(q)
\end{equation}

Then with $(p_1,r_1) \in E(q)$ and $(p_2,r_2) \in E(q)$ \\
we have the group law \\

\begin{equation}
(p_1,r_1)(p_2,r_2) = ( \frac{1}{q} \frac{p_1r_2+p_2r_1}{1-p_1p_2r_1r_2}, \frac{1}{q} \frac{p_1p_2 + r_1r_2}{1+p_1p_2r_1r_2} )
\end{equation}

Then inverse is given by \\
\begin{equation}
(p,r)^{-1} = (-p,r)
\end{equation}

and the identity is \\
\begin{equation}
(0,q)
\end{equation}

\section{Example}
\begin{align*}
& q=\frac{16}{21} \\
& p=\frac{4}{13}, \hspace{0.25cm} r=\frac{4}{5}
\end{align*}

Then $(p,r) \in E(q)$ is a rational point. \bigskip \\
With $\bar{p}=\frac{9}{17}$ and $(\bar{p},s) \in E(q)$, we find $s = \frac{\sqrt{474993}}{801}$ \bigskip \\
Thus $(\bar{p},s)$ is not a rational point on $E(q)$

\bigskip

\noindent\textit{Department of Physics, University of Colorado Boulder, Boulder, CO 80309\\
Walter.Wyss@Colorado.EDU}

\end{document}